\documentclass{article}
\usepackage{amssymb}
\usepackage{amsfonts}
\usepackage{amsmath}

\input{tcilatex}

\begin{document}

\begin{center}
{\large Some remarks regarding Quaternions and Octonions }

\begin{equation*}
\end{equation*}%
Cristina FLAUT\ 

\bigskip 
\begin{equation*}
\end{equation*}
\end{center}

\textbf{Abstract}. {\small In this paper, we present some applications of
quaternions and octonions. We present the real matrix representations for
complex octonions and some of their properties which can be used in
computations where these elements are involved. Moreover, we give a set of
invertible elements in split quaternion algebras and in split octonion
algebras.}

\begin{center}
\begin{equation*}
\end{equation*}
\end{center}


\medskip

\bigskip 
\begin{equation*}
\end{equation*}

\textbf{1. Introduction}%
\begin{equation*}
\end{equation*}

The mathematical objects have, directly or indirectly, many applications in
our lives. It is very interesting to remark how abstract notions and
abstract theories can have such an influence, leading to a significant
development in various domains.

As examples in this direction, Quaternions and Octonions are some of such
objects which influenced, in a good sense, our lives over the time.

Quaternions were discovered in October 1843 by Sir William Rowan Hamilton,
when he introduced an algebraic system formed by one real part and three
imaginary parts. This system is nothing else than the well known real
quaternion division algebra, ( Hanson 2005, p. 5).

Since quaternions generalize complex numbers, after quaternions were
discovered, a question arised: if this structure can be generalized. In
December 1843, John T. Graves, Hamilton's friend, generalized quaternions to
octonions and gave them the name \textquotedblright
octaves\textquotedblright , obtaining an 8-dimensional algebra.
Independently by Graves, Arthur Cayley \ discovered octonions and published
his result in 1845. Octonions are also called \textit{Cayley numbers, }(
Hanson 2005, p. 9).

A generalized real\ quaternion algebra, $\mathbb{H}\left( \beta _{1},\beta
_{2}\right) $, is an algebra with the elements of the form $%
a=a_{0}+a_{1}e_{1}+a_{2}e_{2}+a_{3}e_{3}$, where $a_{i}\in \mathbb{R},i\in
\{0,1,2,3\}$, and the basis $\{1,e_{1},e_{2},e_{3}\}$, with the
multiplication given in the following table:\medskip

\begin{center}
\begin{tabular}{ccccc}
$\cdot $ & $1$ & $e_{1}$ & $e_{2}$ & $e_{3}$ \\ \hline
$1$ & $1$ & $e_{1}$ & $e_{2}$ & $e_{3}$ \\ 
$e_{1}$ & $e_{1}$ & $-\beta _{1}$ & $e_{3}$ & $-\beta _{1}e_{2}$ \\ 
$e_{2}$ & $e_{2}$ & $-e_{3}$ & $-\beta _{2}$ & $\beta _{2}e_{1}$ \\ 
$e_{3}$ & $e_{3}$ & $\beta _{1}e_{2}$ & $-\beta _{2}e_{1}$ & $-\beta
_{1}\beta _{2}$%
\end{tabular}

{\small Table 1.}\medskip
\end{center}

If $a\in \mathbb{H}\left( \beta _{1},\beta _{2}\right) $, $%
a=a_{0}+a_{1}e_{1}+a_{2}e_{2}+a_{3}e_{3}$, then $\bar{a}%
=a_{0}-a_{1}e_{1}-a_{2}e_{2}-a_{3}e_{3}$ is called the \textit{conjugate} of
the element $a$. We denote by 
\begin{equation*}
\boldsymbol{t}\left( a\right) =a+\overline{a}\in \mathbb{R}
\end{equation*}%
and 
\begin{equation*}
\boldsymbol{n}\left( a\right) =a\overline{a}\in \mathbb{R},
\end{equation*}%
\textit{the} \textit{trace} and \textit{the norm} of a real quaternion $a.$
The norm of a generalized quaternion has the following expression 
\begin{equation*}
\boldsymbol{n}\left( a\right) =a_{1}^{2}+\beta _{1}a_{2}^{2}+\beta
_{2}a_{3}^{2}+\beta _{1}\beta _{2}a_{4}^{2}.
\end{equation*}%
\qquad\ \ \ \ \ \ \ \ 

A generalized octonion algebra $\mathbb{O}(\alpha ,\beta ,\gamma )$ is an
algebra with the elements of the form$%
~a=a_{0}+a_{1}e_{1}+a_{2}e_{2}+a_{3}e_{3}+a_{4}e_{4}+a_{5}e_{5}+a_{6}e_{6}+a_{7}e_{7}, 
$ where $a_{i}\in \mathbb{R},i\in \{0,1,2,3,4,5,6,7\},$ and the basis $%
\{1,e_{1},...,e_{7}\},\medskip $ with multiplication given in the following
table:\medskip

\begin{center}
{\footnotesize $%
\begin{tabular}{l|llllllll}
$\cdot $ & $1$ & $\,\,\,e_{1}$ & $\,\,\,\,\,e_{2}$ & $\,\,\,\,e_{3}$ & $%
\,\,\,\,e_{4}$ & $\,\,\,\,\,\,e_{5}$ & $\,\,\,\,\,\,e_{6}$ & $%
\,\,\,\,\,\,\,e_{7}$ \\ \hline
$\,1$ & $1$ & $\,\,\,e_{1}$ & $\,\,\,\,e_{2}$ & $\,\,\,\,e_{3}$ & $%
\,\,\,\,e_{4}$ & $\,\,\,\,\,\,e_{5}$ & $\,\,\,\,\,e_{6}$ & $%
\,\,\,\,\,\,\,e_{7}$ \\ 
$\,e_{1}$ & $\,\,e_{1}$ & $-\alpha $ & $\,\,\,\,e_{3}$ & $-\alpha e_{2}$ & $%
\,\,\,\,e_{5}$ & $-\alpha e_{4}$ & $-\,\,e_{7}$ & $\,\,\,\alpha e_{6}$ \\ 
$\,e_{2}$ & $\,e_{2}$ & $-e_{3}$ & $-\,\beta $ & $\,\,\beta e_{1}$ & $%
\,\,\,\,e_{6}$ & $\,\,\,\,\,e_{7}$ & $-\beta e_{4}$ & $-\beta e_{5}$ \\ 
$e_{3}$ & $e_{3}$ & $\alpha e_{2}$ & $-\beta e_{1}$ & $-\alpha \beta $ & $%
\,\,\,\,e_{7}$ & $-\alpha e_{6}$ & $\,\,\,\beta e_{5}$ & $-\alpha \beta
e_{4} $ \\ 
$e_{4}$ & $e_{4}$ & $-e_{5}$ & $-\,e_{6}$ & $-\,\,e_{7}$ & $-\,\gamma $ & $%
\,\,\,\gamma e_{1}$ & $\,\,\gamma e_{2}$ & $\,\,\,\,\,\gamma e_{3}$ \\ 
$\,e_{5}$ & $\,e_{5}$ & $\alpha e_{4}$ & $-\,e_{7}$ & $\,\alpha e_{6}$ & $%
-\gamma e_{1}$ & $-\,\alpha \gamma $ & $-\gamma e_{3}$ & $\,\alpha \gamma
e_{2}$ \\ 
$\,\,e_{6}$ & $\,\,e_{6}$ & $\,\,\,\,e_{7}$ & $\,\,\beta e_{4}$ & $-\,\beta
e_{5}$ & $-\gamma e_{2}$ & $\,\,\,\gamma e_{3}$ & $-\beta \gamma $ & $-\beta
\gamma e_{1}$ \\ 
$\,\,e_{7}$ & $\,\,e_{7}$ & $-\alpha e_{6}$ & $\,\beta e_{5}$ & $\alpha
\beta e_{4}$ & $-\gamma e_{3}$ & $-\alpha \gamma e_{2}$ & $\beta \gamma
e_{1} $ & $-\alpha \beta \gamma $%
\end{tabular}%
\ \medskip $ }

{\small Table 2.}\medskip
\end{center}

The algebra $\mathbb{O}(\alpha ,\beta ,\gamma )$ is non-commutative,
non-associative but it is \textit{alternative}\ (i.e.\thinspace \thinspace $%
x^{2}y=x\left( xy\right) $ and $yx^{2}=\left( yx\right) x,\forall x,y\in 
\mathbb{O}(\alpha ,\beta ,\gamma )$), \textit{flexible}\ (i.e. $x\left(
yx\right) =\left( xy\right) x,\forall x,y\in \mathbb{O}(\alpha ,\beta
,\gamma )$) and \textit{power-associative}\ (i.e. for each $x\in \mathbb{O}%
(\alpha ,\beta ,\gamma )$, the subalgebra generated by $x$ is an associative
algebra). For other details regarding quaternions and octonions, the reader
is referred to (Schafer, 1966).

If $a\in \mathbb{O}(\alpha ,\beta ,\gamma ),$ $%
a=a_{0}+a_{1}e_{1}+a_{2}e_{2}+a_{3}e_{3}+a_{4}e_{4}+a_{5}e_{5}+a_{6}e_{6}+a_{7}e_{7}, 
$ then $\bar{a}%
=a_{0}-a_{1}e_{1}-a_{2}e_{2}-a_{3}e_{3}-a_{4}e_{4}-a_{5}e_{5}-a_{6}e_{6}-a_{7}e_{7} 
$ is called the \textit{conjugate} of the element $a.$ Let $a\in \mathbb{O}%
(\alpha ,\beta ,\gamma ).$ We have the \textit{trace}, respectively, the 
\textit{norm} of the element $a$ given by the relations 
\begin{equation*}
\mathbf{t}\left( a\right) =a+\overline{a}\in \mathbb{R}
\end{equation*}%
and 
\begin{equation*}
\,\mathbf{n}\left( a\right) =a\overline{a}=a_{0}^{2}+\alpha a_{1}^{2}+\beta
a_{2}^{2}+\alpha \beta a_{3}^{2}+\gamma a_{4}^{2}+\alpha \gamma
a_{5}^{2}+\beta \gamma a_{6}^{2}+\alpha \beta \gamma a_{7}^{2}\in \mathbb{R}.
\end{equation*}

We remark that the following relation holds $\,$\thinspace \thinspace 
\begin{equation*}
a^{2}-\mathbf{t}\left( a\right) a+\mathbf{n}\left( a\right) =0,
\end{equation*}%
for each $a\in A,~$where $A\in \{\mathbb{H}\left( \beta _{1},\beta
_{2}\right) ,\mathbb{O}(\alpha ,\beta ,\gamma )\}.$

We know that a finite-dimensional algebra $A$ is \textit{a division} algebra
if and only if $A$ does not contain zero divisors, (Schafer, 1966). With the
above notations, for $\beta _{1}=\beta _{2}=1,$ we obtain the real division
algebra $\mathbb{H}$ and for $\alpha =\beta =\gamma =1,$ we obtain the real
division octonion algebra $\mathbb{O}.~$

If a quaternion algebra and an octonion algebra are not division algebras,
we call them \textit{a split} \textit{quaternion algebra}, respectively, 
\textit{a split octonion algebra}.

From the above properties, we can see that an algebra $A,~$with \newline
$A\in \{\mathbb{H}\left( \beta _{1},\beta _{2}\right) ,\mathbb{O}(\alpha
,\beta ,\gamma )\},$ is a division algebra if and only if $\mathbf{n}\left(
x\right) \neq 0,$ for each $x\in A,x\neq 0$, (Schafer, 1966, p. 27).
Therefore, depending on the values of the numbers $\beta _{1},\beta
_{2},\alpha ,\beta ,\gamma ,$ we can obtain division algebras or split
algebras.

A real quaternion algebra is isomorphic to $\mathbb{H}$, when we have a
division algebra, or it \ is isomorphic to the algebra $\mathcal{M}%
_{2}\left( R\right) ,$ the algebra of~ $2\times 2$ real matrices, when we
have a split algebra, (Schafer, 1966, p. 25).

We consider the following real algebra 
\begin{equation*}
\mathbb{A}=\{\left( 
\begin{array}{ll}
a & u \\ 
v & b%
\end{array}%
\right) ,a,b\in \mathbb{R},\,u,v\in \mathbb{R}^{3}\},
\end{equation*}%
with the usual addition and scalar multiplication of matrices. We consider
the following multiplication 
\begin{equation}
\left( 
\begin{array}{ll}
a & u \\ 
v & b%
\end{array}%
\right) \left( 
\begin{array}{ll}
c & z \\ 
w & d%
\end{array}%
\right) =\left( 
\begin{array}{ll}
\,\,\,\,ac+\left( u,w\right) & az+du-v\times w \\ 
cv+bw+u\times z & \,\,\,\,\,\,\,\,\,\,bd+<v,z>%
\end{array}%
\right) ,  \tag{1.1}
\end{equation}%
with $<,>$ the usual dot product and $\times $ the cross product of vectors
from $\mathbb{R}^{3}.$ With matrices addition, scalar multiplication of
matrices and the multiplication given in relation $\left( 1.1\right) ,$ we
obtain a unitary non-associative algebra of dimension $8$, called \ the 
\textit{Zorn's vector-matrix algebra}. Therefore, a real octonion algebra is
isomorphic to $\mathbb{O}$, when we have a division algebra, or it \ is
isomorphic to the algebra $\mathbb{A}$, when we have a split algebra, see
(Kostrikin, Shafarevich, 1995). A famous Hurwitz's\textbf{\ }theorem states
that $\mathbb{R}$, $\mathbb{C}$, $\mathbb{H}$ and $\mathbb{O}$ are the only
real alternative division algebras, see (Baez, 2002).

Division algebras are used to built space-time block codes and division
quaternion algebras are also used for this purpose. Quaternions are used in
Coding Theory and in digital signal processing, see (Alfsmann et all, 2007).
There are a lot of examples of codes based on quaternion algebras, first of
them, Alamouti code, appeared in 1998, see (Alamouti, 1998) and (Unger,
Markin, 2011).

For other conections of different types of algebras with codes, the reader
is referred to (Flaut, 2015).

Another application of quaternions is in representing rotations in $\mathbb{R%
}^{3}.\,$\ Usually, a rotation in $\mathbb{R}^{3}$ around an axis is given
by a square orthogonal matrix of order three, with its determinant equal
with $1.$ When we compute two rotations, a lot of computations are involved.
\ In this situations, quaternions give a much better representation for a
rotation. In data registration, when we want to find a transformation for a
set of data points such that these points fit better to a shape model,
quaternions are used for solving this problem, since it is necessary to find
a rotation and a translation with some good properties, see (Jia, 2017).

Octonions have many applications in processing of color images, as for
example in color image edge detection, see (Chen, Tu, 2014), in remote
sensing images, (Li, 2011), in 2D and 3D signal processing, (Snopek, 15), in
artificial neural networks, where Octonionic neural networks are used as a
computational models, (Klco et all, 2017),\ \ in electrodynamics, see
(Chanyal, 2011) and (Chen, Tu, 2014), in wireless data communication
(Jouget, 2013) and examples can continue.

Since the above applications used quaternion and octonion matrices or matrix
representations of these algebras, in this chapter, we present real matrix
representations for complex octonions and \ some of their properties and we
provide examples of invertible elements in split quaternion algebras and in
split octonion algebras.

For other details regarding quaternions and octonions, the reader are
referred to (Schafer, 1966, p. 27), (Flaut, Savin, 2015), (Flaut, Savin,
2014), (Flaut, Savin, Iorgulescu 2013), (Flaut, \c{S}tefanescu, 2009),
(Flaut, 2006), (Flaut, Shpakivskyi, 2013(1)), (Flaut, Shpakivskyi, 2015),
(Flaut, Shpakivskyi, 2015(1)), (Flaut, Shpakivskyi, 2015(1)). 
\begin{equation*}
\end{equation*}%
2. \textbf{Real matrix representations for complex octonions}

\begin{equation*}
\end{equation*}%
\qquad\ 

We know that in an alternative algebra $A$, the following identities, called 
\textit{Moufang identities}, are true:%
\begin{equation}
\left( xzx\right) y=x[z\left( xy\right) ],  \tag{2.1.}
\end{equation}%
\begin{equation}
y\left( xzx\right) =[\left( yx\right) z]x,  \tag{2.2.}
\end{equation}%
\begin{equation}
\left( xy\right) \left( zx\right) =x\left( yz\right) x,  \tag{2.3.}
\end{equation}%
for all $x,y,z\in A.$ Since octonions \ form an alternative algebra, we can
use these relations in computations.

Let $\mathbb{O}$ be the real\ division octonion algebra, the algebra of the
elements of the form $%
a=a_{0}+a_{1}i+a_{2}j+a_{3}ij+a_{4}k+a_{5}ik+a_{6}jk+a_{7}(ij)k,$ where

\begin{equation*}
a_{i}\in \mathbb{R},i^{2}=j^{2}=k^{2}=-1,
\end{equation*}%
and \ 
\begin{equation*}
ij=-ji,jk=-kj,ki=-ik,\left( ij\right) k=-i(jk).
\end{equation*}
The set $\{1,i,j,ij,k,ik,jk,\left( ij\right) k\}$ is a basis in $\mathbb{O}$.

We consider $K$ to be the field $\{\left( 
\begin{array}{cc}
a & -b \\ 
b & a%
\end{array}%
\right) $ $\mid $ $a,b\in \mathbb{R}\}$ and the map 
\begin{equation*}
\varphi :\mathbb{C}\rightarrow K,\varphi \left( a+bi\right) =\left( 
\begin{array}{cc}
a & -b \\ 
b & a%
\end{array}%
\right) ,
\end{equation*}%
where $i^{2}=-1.$ The map $\varphi $ is a fields morphism and the element $%
\varphi \left( z\right) =\left( 
\begin{array}{cc}
a & -b \\ 
b & a%
\end{array}%
\right) $ is called the matrix representation of the complex element $%
z=a+bi\in \mathbb{C}.$

A complex octonion is an element of the form $%
A=a_{0}+a_{1}e_{1}+a_{2}e_{2}+a_{3}e_{3}+a_{4}e_{4}+a_{5}e_{5}+a_{6}e_{6}+a_{7}e_{7}, 
$ where $a_{m}\in \mathbb{C},m\in \{0,1,2,...,7\},$%
\begin{equation*}
\ e_{m}^{2}=-1,\,\,\,m\in \{0,1,2,...,7\}
\end{equation*}%
and \ 
\begin{equation*}
e_{m}e_{n}=-e_{n}e_{m}=\gamma _{mn}e_{t},\,\,\gamma _{mn}\in \{-1,1\},m\neq
n,m,n\in \{0,1,2,...,7\},
\end{equation*}%
$\gamma _{mn}$ and $e_{t}$ are uniquely determined by $e_{m}$ and $e_{n}.$
We denote by $\mathbb{O}_{C}$ the algebra of the complex octonions\textit{,}
called \textit{the complex octonion algebra}. This algebra is an algebra
over the field $\mathbb{C}$ and the set $%
\{1,e_{1},e_{2},e_{3},e_{4},e_{5},e_{6},e_{7}\}$ is a basis in $\mathbb{O}%
_{C}$.

The map $\delta :\mathbb{R}\rightarrow \mathbb{C},\delta \left( a\right) =a$
is the inclusion morphism between $\mathbb{R}$ and $\mathbb{C}$ as $\mathbb{R%
}$-algebras$.$ We denote by $\mathbb{G}$ the following $\mathbb{C}$%
-subalgebra of the algebra $\mathbb{O}_{\mathbb{C}},$%
\begin{equation*}
\mathbb{G}=\{A\in \mathbb{O}_{C}\ \mid ~A=\overset{7}{\underset{m=0}{\sum }}%
a_{m}e_{m},a_{m}\in \mathbb{R},m\in \{0,1,2,3,...,7\}\}.
\end{equation*}

It is clear that $\mathbb{G}$ becomes an algebra over $\mathbb{R},$ with the
following multiplication $"\cdot "$ 
\begin{equation*}
a\cdot A=\delta \left( a\right) A=aA,a\in \mathbb{R},A\in \mathbb{G}.
\end{equation*}%
We denote this algebra by $\mathbb{O}_{\mathbb{R}}$ and the map 
\begin{eqnarray*}
\phi &:&\mathbb{O\rightarrow O}_{\mathbb{R}},\phi \left(
a_{0}+a_{1}i+a_{2}j+a_{3}ij+a_{4}k+a_{5}ik+a_{6}jk+a_{7}(ij)k\right) = \\
&=&a_{0}+a_{1}e_{1}+a_{2}e_{2}+a_{3}e_{3}+a_{4}e_{4}+a_{5}e_{5}+a_{6}e_{6}+a_{7}e_{7},
\end{eqnarray*}%
where $a_{m}\in \mathbb{R},m\in \{0,1,2,3,...,7\}$ is an algebra
isomorphism. Due to this isomorphism, we have that $\phi \left( 1\right)
=1,\phi \left( i\right) =e_{1},\phi \left( j\right) =e_{2},\phi \left(
ij\right) =e_{3},$\newline
$\phi \left( k\right) =e_{4},\phi \left( ik\right) =e_{5},\phi \left(
jk\right) =e_{6},\phi \left( \left( ij\right) k\right) =e_{7}$ and the
algebras $\mathbb{O}_{\mathbb{R}}$ and $\mathbb{O}_{C}$ have the same basis $%
\{1,e_{1},e_{2},e_{3},...,e_{7}\}.$ With these notations, in the rest of the
paper, we denote the real octonion $%
a_{0}+a_{1}i+a_{2}j+a_{3}ij+a_{4}k+a_{5}ik+a_{6}jk+a_{7}(ij)k$~~with the
octonion $%
a_{0}+a_{1}e_{1}+a_{2}e_{2}+a_{3}e_{3}+a_{4}e_{4}+a_{5}e_{5}+a_{6}e_{6}+a_{7}e_{7} 
$, and viceversa, where $a_{m}\in \mathbb{R},m\in \{0,1,2,3,...,7\}$ and we
use the notation $\mathbb{O}$ instead of $\mathbb{O}_{\mathbb{R}}.$

We consider the complex octonion\newline
$%
A=a_{0}+a_{1}e_{1}+a_{2}e_{2}+a_{3}e_{3}+a_{4}e_{4}+a_{5}e_{5}+a_{6}e_{6}+a_{7}e_{7}, 
$\newline
$a_{m}\in \mathbb{C},m\in \{0,1,2,3,...,7\}.$ This octonion can be written
under the following form \newline
$%
A=(x_{0}+iy_{0})+(x_{1}+iy_{1})e_{1}+(x_{2}+iy_{2})e_{2}+(x_{3}+iy_{3})e_{3}+ 
$\newline
$%
+(x_{4}+iy_{4})e_{4}+(x_{5}+iy_{5})e_{5}+(x_{6}+iy_{6})e_{6}+(x_{7}+iy_{7})e_{7}, 
$ where $x_{m},y_{m}\in \mathbb{R},$\newline
$m\in \{0,1,2,3,...,7\}$ and $i^{2}=-1.$

It results that 
\begin{equation}
A=x+iy,  \tag{2.4.}
\end{equation}%
with $x,y\in \mathbb{O}$, $%
x=x_{0}+x_{1}e_{1}+x_{2}e_{2}+x_{3}e_{3}+x_{4}e_{4}+x_{5}e_{5}+x_{6}e_{6}+x_{7}e_{7} 
$ and \newline
$%
y=y_{0}+y_{1}e_{1}+y_{2}e_{2}+y_{3}e_{3}+y_{4}e_{4}+y_{5}e_{5}+y_{6}e_{6}+y_{7}e_{7}. 
$

The conjugate of this octonion is the element $\overline{A}%
=a_{0}-a_{1}e_{1}-a_{2}e_{2}-a_{3}e_{3}-a_{4}e_{4}-a_{5}e_{5}-a_{6}e_{6}-a_{7}e_{7}, 
$ therefore, with the above notations, we have 
\begin{equation}
\overline{A}=\overline{x}+i\overline{y}.  \tag{2.5.}
\end{equation}

To an octonion $%
a=a_{0}+a_{1}e_{1}+a_{2}e_{2}+a_{3}e_{3}+a_{4}e_{4}+a_{5}e_{5}+a_{6}e_{6}+a_{7}e_{7} 
$ $\in \mathbb{O}$, we associate the following element%
\begin{equation}
a^{\ast
}=a_{0}+a_{1}e_{1}-a_{2}e_{2}-a_{3}e_{3}-a_{4}e_{4}-a_{5}e_{5}-a_{6}e_{6}-a_{7}e_{7}.
\tag{2.6.}
\end{equation}%
We remark that 
\begin{equation}
(a^{\ast })^{\ast }=a  \tag{2.7}
\end{equation}%
and%
\begin{equation}
\left( a+b\right) ^{\ast }=a^{\ast }+b^{\ast },  \tag{2.8}
\end{equation}

The identity $(ab)^{\ast }=a^{\ast }b^{\ast },~$for all $a,b\in \mathbb{O,}$
in general is not true. Since the real octonion $%
a_{0}+a_{1}i+a_{2}j+a_{3}ij+a_{4}k+a_{5}ik+a_{6}jk+a_{7}(ij)k$ can be
written under the form $a=q_{1}+q_{2}k,$ where $%
q_{1}=a_{0}+a_{1}i+a_{2}j+a_{3}ij$ and $q_{2}=a_{4}+a_{5}ei+a_{6}j+a_{7}ij$
are two real quaternions, it results 
\begin{equation}
a^{\ast }=q_{1}^{\ast }-q_{2}k,  \tag{2.9}
\end{equation}%
where $q_{1}^{\ast }=a_{0}+a_{1}i-a_{2}j-a_{3}ij.$

With the above notations, we define

\begin{equation}
\widetilde{a}=q_{1}-q_{2}k,  \tag{2.10}
\end{equation}%
\begin{equation}
a_{+}=q_{1}^{\ast }+q_{2}^{\ast }k,~a^{+}=\overline{q}_{1}+q_{2}k. 
\tag{2.11}
\end{equation}

Considering \ the real division quaternion algebra $\mathbb{H},\,$in (Tian,
2000) were defined $\ $the following maps \ 
\begin{equation}
\lambda :\mathbb{H}\rightarrow \mathcal{M}_{4}\left( \mathbb{R}\right)
,\lambda \left( a\right) =\left( 
\begin{array}{llll}
a_{0} & -a_{1} & -a_{2} & -a_{3} \\ 
a_{1} & a_{0} & -a_{3} & a_{2} \\ 
a_{2} & a_{3} & a_{0} & -a_{1} \\ 
a_{3} & -a_{2} & a_{1} & a_{0}%
\end{array}%
\right)  \tag{2.12}
\end{equation}%
and 
\begin{equation}
\rho :\mathbb{H}\rightarrow \mathcal{M}_{4}\left( \mathbb{R}\right)
,\,\,\,\rho \left( a\right) =\left( 
\begin{array}{llll}
a_{0} & -a_{1} & -a_{2} & -a_{3} \\ 
a_{1} & a_{0} & a_{3} & -a_{2} \\ 
a_{2} & -a_{3} & a_{0} & a_{1} \\ 
a_{3} & a_{2} & -a_{1} & a_{0}%
\end{array}%
\right) ,\medskip  \tag{2.13}
\end{equation}%
where $a=a_{0}+a_{1}i+a_{2}j+a_{3}ij\in \mathbb{H.}\,\,$\ 

We remark that $\lambda $ is an isomorphism between $\mathbb{H}\,\ $and the
algebra of the \ matrices 
\begin{equation*}
\left\{ \left( 
\begin{array}{llll}
a_{0} & -a_{1} & -a_{2} & -a_{3} \\ 
a_{1} & a_{0} & -a_{3} & a_{2} \\ 
a_{2} & a_{3} & a_{0} & -a_{1} \\ 
a_{3} & -a_{2} & a_{1} & a_{0}%
\end{array}%
\right) ,a_{0},a_{1},a_{2},a_{3}\in \mathbb{R}\right\} .\,\,
\end{equation*}%
The columns of the matrix $\lambda \left( a\right) \in \mathcal{M}_{4}\left( 
\mathbb{R}\right) $ are represented by the coefficients in $\mathbb{R}$ of
the elements $\ \{a,ai,aj,a\left( ij\right) \},$\ considered in respect to
the basis $\{1,i,j,ij\}.$

The matrix $\lambda \left( a\right) $ is called \textit{the left matrix
representation} of the element $a\in \mathbb{H}.\medskip $

The map $\rho $ is an isomorphism between $\mathbb{H}\,\ $and the algebra of
the \ matrices 
\begin{equation*}
\left\{ \left( 
\begin{array}{llll}
a_{0} & -a_{1} & -a_{2} & -a_{3} \\ 
a_{1} & a_{0} & a_{3} & -a_{2} \\ 
a_{2} & -a_{3} & a_{0} & a_{1} \\ 
a_{3} & a_{2} & -a_{1} & a_{0}%
\end{array}%
\right) ,a_{0},a_{1},a_{2},a_{3}\in \mathbb{R}\right\} .\,\,
\end{equation*}

In a similar way, we remark that the columns of the matrix $\rho \left(
a\right) \in \mathcal{M}_{4}\left( \mathbb{R}\right) $ are the coefficients
in $\mathbb{R}$ of the elements $\{a,ia,ja,\left( ij\right) a\},$ considered
in respect to the basis $\{1,i,j,ij\}.$ The matrix $\rho \left( a\right) $
is called \textit{the right matrix representation} of the quaternion $a\in 
\mathbb{H}.\medskip $

With these notations, in (Tian, 2000) were defined the left and right real
representations of the octonion $a=q_{1}+q_{2}k$, namely 
\begin{equation}
\Lambda \left( a\right) =\left( 
\begin{array}{cc}
\lambda \left( q_{1}\right) & -\rho \left( q_{2}\right) M_{1} \\ 
\lambda \left( q_{2}\right) M_{1} & \rho \left( q_{1}\right)%
\end{array}%
\right) \in \mathcal{M}_{8}\left( \mathbb{R}\right)  \tag{2.14}
\end{equation}

and 
\begin{equation}
\Delta \left( a\right) =\left( 
\begin{array}{cc}
\rho \left( q_{1}\right) & -\lambda \left( \overline{q}_{2}\right) \\ 
\lambda \left( q_{2}\right) & \rho \left( \overline{q}_{1}\right)%
\end{array}%
\right) \in \mathcal{M}_{8}\left( \mathbb{R}\right) ,  \tag{2.15}
\end{equation}%
where $%
q_{1}=a_{0}+a_{1}i+a_{2}j+a_{3}ij,~q_{2}=a_{4}+a_{5}ei+a_{6}j+a_{7}ij~ $ are
two real quaternions and $M_{1}=\left( 
\begin{array}{cccc}
1 & 0 & 0 & 0 \\ 
0 & -1 & 0 & 0 \\ 
0 & 0 & -1 & 0 \\ 
0 & 0 & 0 & -1%
\end{array}%
\right) $.~We remark that $L_{1}=\lambda \left( i\right) =\left( 
\begin{array}{llll}
0 & -1 & 0 & 0 \\ 
1 & 0 & 0 & 0 \\ 
0 & 0 & 0 & -1 \\ 
0 & 0 & 1 & 0%
\end{array}%
\right) $ and $R_{1}=\rho \left( i\right) =\left( 
\begin{array}{llll}
0 & -1 & 0 & 0 \\ 
1 & 0 & 0 & 0 \\ 
0 & 0 & 0 & 1 \\ 
0 & 0 & -1 & 0%
\end{array}%
\right) .$

In (Flaut, Shpakivskyi, 2013), using the above notations, were obtained the
following matrices 
\begin{equation}
\Gamma \left( Q\right) =\left( 
\begin{array}{cc}
\lambda \left( a\right) & -\lambda \left( b^{\ast }\right) \\ 
\lambda \left( b\right) & \lambda \left( a^{\ast }\right)%
\end{array}%
\right)  \tag{2.16}
\end{equation}%
and 
\begin{equation}
\Theta \left( Q\right) =\left( 
\begin{array}{cc}
\rho \left( a\right) & -\rho \left( b\right) \\ 
\rho \left( b^{\ast }\right) & \rho \left( a^{\ast }\right)%
\end{array}%
\right) ,  \tag{2.17}
\end{equation}%
where $Q=a+ib$ is a complex quaternion, with $%
a=a_{0}+a_{1}e_{1}+a_{2}e_{2}+a_{3}e_{3}\in \mathbb{H}%
,b=b_{0}+b_{1}e_{1}+b_{2}e_{2}+b_{3}e_{3}\in \mathbb{H,}$ $a^{\ast
}=a_{0}+a_{1}e_{1}-a_{2}e_{2}-a_{3}e_{3}\in \mathbb{H,~}b^{\ast
}=b_{0}+b_{1}e_{1}-b_{2}e_{2}-b_{3}e_{3}\in \mathbb{H~~}$and $i^{2}=-1.$ The
matrix $\Gamma \left( Q\right) \in \mathcal{M}_{8}\left( \mathbb{R}\right) $
is called \textit{the left real matrix representation for the complex
quaternion} $Q$ and $\Theta \left( Q\right) \in \mathcal{M}_{8}\left( 
\mathbb{R}\right) $ is called the \textit{right real matrix representation}
for the \ complex quaternion $Q$.

Using some ideas developed above, in the following, we define the left and
the right real matrix representations for the complex octonions\ and we
investigate some of their properties.

We consider the complex octonion $A=x+iy,$with $x,y\in \mathbb{O}$, $%
x=x_{0}+x_{1}e_{1}+x_{2}e_{2}+x_{3}e_{3}+x_{4}e_{4}+x_{5}e_{5}+x_{6}e_{6}+x_{7}e_{7} 
$ and \newline
$%
y=y_{0}+y_{1}e_{1}+y_{2}e_{2}+y_{3}e_{3}+y_{4}e_{4}+y_{5}e_{5}+y_{6}e_{6}+y_{7}e_{7}. 
$

The matrices%
\begin{equation}
\Phi \left( A\right) =\left( 
\begin{array}{cc}
\Lambda \left( x\right) & -\Lambda \left( y\right) \\ 
\Lambda \left( y^{\ast }\right) & \Lambda \left( x^{\ast }\right)%
\end{array}%
\right) \in \mathcal{M}_{16}\left( \mathbb{R}\right) ,  \tag{2.18}
\end{equation}

$~$and 
\begin{equation}
\Psi \left( A\right) =\left( 
\begin{array}{cc}
\Delta \left( x\right) & -\Delta \left( y\right) \\ 
\Delta \left( y^{\ast }\right) & \Delta \left( x^{\ast }\right)%
\end{array}%
\right) \in \mathcal{M}_{16}\left( \mathbb{R}\right)  \tag{2.19}
\end{equation}%
are called \textit{the left real matrix representation} and \textit{the
right real matrix representation }for the complex octonion $A.\medskip $

\textbf{Definition 2.1}. Let $A\in \mathbb{O}_{\mathbb{C}},A=x+iy.$ We
consider the following matrix 
\begin{equation*}
\overrightarrow{A}=(\overrightarrow{x}^{t},\overrightarrow{y}%
^{t})^{t}=\left( 
\begin{array}{c}
\overrightarrow{x} \\ 
\overrightarrow{y}%
\end{array}%
\right) \in \mathcal{M}_{16\times 1}\left( \mathbb{R}\right) ,
\end{equation*}%
$~$\newline
called \textit{the vector representation }of the element $X$, where $x,y\in 
\mathbb{O}$, \newline
$%
x=x_{0}+x_{1}e_{1}+x_{2}e_{2}+x_{3}e_{3}+x_{4}e_{4}+x_{5}e_{5}+x_{6}e_{6}+x_{7}e_{7}
$,\newline
$%
y=y_{0}+y_{1}e_{1}+y_{2}e_{2}+y_{3}e_{3}+y_{4}e_{4}+y_{5}e_{5}+y_{6}e_{6}+y_{7}e_{7}
$,\newline
$\overrightarrow{x}=(x_{0},x_{1},x_{2},x_{3},x_{4},x_{5},x_{6},x_{7})^{t}\in 
\mathcal{M}_{8\times 1}\left( \mathbb{R}\right) ,$\newline
$\overrightarrow{y}=(y_{0},y_{1},y_{2},y_{3},y_{4},y_{5},y_{6},y_{7})^{t}\in 
\mathcal{M}_{8\times 1}\left( \mathbb{R}\right) $ are the vector
representations for the real octonions $x$ and $y,$ as were defined in
(Tian, 2000). If $\overrightarrow{x}%
=(x_{0},x_{1},x_{2},x_{3},x_{4},x_{5},x_{6},x_{7})^{t}\in \mathcal{M}%
_{8\times 1}\left( \mathbb{R}\right) $ and\newline
$\overrightarrow{y}=(y_{0},y_{1},y_{2},y_{3},y_{4},y_{5},y_{6},y_{7})^{t}\in 
\mathcal{M}_{8\times 1}\left( \mathbb{R}\right) $, we have 
\begin{equation}
\overrightarrow{xy}=\Lambda \left( x\right) \overrightarrow{y}  \tag{2.20}
\end{equation}%
and 
\begin{equation}
\overrightarrow{yx}=\Delta \left( x\right) \overrightarrow{y}.  \tag{2.21}
\end{equation}

\bigskip From the same paper, Theorem 2.1, Theorem 2.3 and Theorem 2.9, we
have that\newline
\begin{equation}
\overrightarrow{ax}=\Lambda \left( a\right) \overrightarrow{x},%
\overrightarrow{xa}=\Delta \left( a\right) \overrightarrow{x},  \tag{2.22}
\end{equation}%
\begin{equation}
\Lambda \left( a^{2}\right) =\Lambda \left( a\right) \Lambda \left( a\right)
,\Delta \left( a^{2}\right) =\Delta \left( a\right) \Delta \left( a\right) ,
\tag{2.23}
\end{equation}%
\begin{equation}
\Lambda \left( a\right) \Delta \left( a\right) =\Delta \left( a\right)
\Lambda \left( a\right) ,  \tag{2.24}
\end{equation}%
where $a,x$ are real octonions.

\textbf{Remark 2.2.} With the above notations, $\,$we have 
\begin{equation*}
\varepsilon \overrightarrow{x}=\overrightarrow{\widetilde{x}}
\end{equation*}%
and 
\begin{equation*}
\tau \overrightarrow{x}=\overrightarrow{x^{\ast }},
\end{equation*}%
where $\varepsilon $=$diag(1,1,1,1,-1,-1,-1,-1)\in \mathcal{M}_{8}\left( 
\mathbb{R}\right) ,$\newline
$\tau $=$diag(1,1,-1,-1,1,1,-1,-1)\in \mathcal{M}_{8}\left( \mathbb{R}%
\right) .\medskip $ Indeed, \newline
$\varepsilon \overrightarrow{x}$=$\left( 
\begin{array}{cccccccc}
1 & 0 & 0 & 0 & 0 & 0 & 0 & 0 \\ 
0 & 1 & 0 & 0 & 0 & 0 & 0 & 0 \\ 
0 & 0 & 1 & 0 & 0 & 0 & 0 & 0 \\ 
0 & 0 & 0 & 1 & 0 & 0 & 0 & 0 \\ 
0 & 0 & 0 & 0 & -1 & 0 & 0 & 0 \\ 
0 & 0 & 0 & 0 & 0 & -1 & 0 & 0 \\ 
0 & 0 & 0 & 0 & 0 & 0 & -1 & 0 \\ 
0 & 0 & 0 & 0 & 0 & 0 & 0 & -1%
\end{array}%
\right) \left( 
\begin{array}{c}
y_{0} \\ 
y_{1} \\ 
y_{2} \\ 
y_{3} \\ 
y_{4} \\ 
y_{5} \\ 
y_{6} \\ 
y_{7}%
\end{array}%
\right) $= $\left( 
\begin{array}{c}
y_{0} \\ 
y_{1} \\ 
y_{2} \\ 
y_{3} \\ 
-y_{4} \\ 
-y_{5} \\ 
-y_{6} \\ 
-y_{7}%
\end{array}%
\right) $=$\overrightarrow{\widetilde{x}}$.\medskip\ In the same way we
obtain\newline
$\tau \overrightarrow{x}$= $\left( 
\begin{array}{cccccccc}
1 & 0 & 0 & 0 & 0 & 0 & 0 & 0 \\ 
0 & 1 & 0 & 0 & 0 & 0 & 0 & 0 \\ 
0 & 0 & -1 & 0 & 0 & 0 & 0 & 0 \\ 
0 & 0 & 0 & -1 & 0 & 0 & 0 & 0 \\ 
0 & 0 & 0 & 0 & 1 & 0 & 0 & 0 \\ 
0 & 0 & 0 & 0 & 0 & 1 & 0 & 0 \\ 
0 & 0 & 0 & 0 & 0 & 0 & -1 & 0 \\ 
0 & 0 & 0 & 0 & 0 & 0 & 0 & -1%
\end{array}%
\right) \left( 
\begin{array}{c}
y_{0} \\ 
y_{1} \\ 
y_{2} \\ 
y_{3} \\ 
y_{4} \\ 
y_{5} \\ 
y_{6} \\ 
y_{7}%
\end{array}%
\right) =\allowbreak \left( 
\begin{array}{c}
y_{0} \\ 
y_{1} \\ 
-y_{2} \\ 
-y_{3} \\ 
y_{4} \\ 
y_{5} \\ 
-y_{6} \\ 
-y_{7}%
\end{array}%
\right) $=$\overrightarrow{x_{+}}.\Box $

\bigskip \textbf{Proposition 2.3.} \textit{Using  the above notations and
definitions, we have}%
\begin{equation*}
\sigma \lambda \left( q_{1}\right) \sigma =\lambda \left( q_{1}^{\ast
}\right) 
\end{equation*}%
\textit{and} 
\begin{equation*}
\sigma \rho \left( q_{1}\right) \sigma =\rho \left( q_{1}^{\ast }\right) ,
\end{equation*}%
\textit{where} $\sigma =diag(1,1,-1,-1).\medskip $

\textbf{Proof.} Indeed, $\sigma \lambda \left( q_{1}\right) \sigma =$\newline
$=$ $\left( 
\begin{array}{cccc}
1 & 0 & 0 & 0 \\ 
0 & 1 & 0 & 0 \\ 
0 & 0 & -1 & 0 \\ 
0 & 0 & 0 & -1%
\end{array}%
\right) \left( 
\begin{array}{llll}
a_{0} & -a_{1} & -a_{2} & -a_{3} \\ 
a_{1} & a_{0} & -a_{3} & a_{2} \\ 
a_{2} & a_{3} & a_{0} & -a_{1} \\ 
a_{3} & -a_{2} & a_{1} & a_{0}%
\end{array}%
\right) \left( 
\begin{array}{cccc}
1 & 0 & 0 & 0 \\ 
0 & 1 & 0 & 0 \\ 
0 & 0 & -1 & 0 \\ 
0 & 0 & 0 & -1%
\end{array}%
\right) =$ $\allowbreak $\newline
$=\left( 
\begin{array}{cccc}
a_{0} & -a_{1} & a_{2} & a_{3} \\ 
a_{1} & a_{0} & a_{3} & -a_{2} \\ 
-a_{2} & -a_{3} & a_{0} & -a_{1} \\ 
-a_{3} & a_{2} & a_{1} & a_{0}%
\end{array}%
\right) \allowbreak =\lambda \left( q_{1}^{\ast }\right) \allowbreak .$%
\newline

In the same way, we obtain $\sigma \rho \left( q_{1}\right) \sigma =$\newline
$=\left( 
\begin{array}{cccc}
1 & 0 & 0 & 0 \\ 
0 & 1 & 0 & 0 \\ 
0 & 0 & -1 & 0 \\ 
0 & 0 & 0 & -1%
\end{array}%
\right) \left( 
\begin{array}{llll}
a_{0} & -a_{1} & -a_{2} & -a_{3} \\ 
a_{1} & a_{0} & a_{3} & -a_{2} \\ 
a_{2} & -a_{3} & a_{0} & a_{1} \\ 
a_{3} & a_{2} & -a_{1} & a_{0}%
\end{array}%
\right) \left( 
\begin{array}{cccc}
1 & 0 & 0 & 0 \\ 
0 & 1 & 0 & 0 \\ 
0 & 0 & -1 & 0 \\ 
0 & 0 & 0 & -1%
\end{array}%
\right) =$\newline
$=\left( 
\begin{array}{cccc}
a_{0} & -a_{1} & a_{2} & a_{3} \\ 
a_{1} & a_{0} & -a_{3} & a_{2} \\ 
-a_{2} & a_{3} & a_{0} & a_{1} \\ 
-a_{3} & -a_{2} & -a_{1} & a_{0}%
\end{array}%
\right) \allowbreak =\rho \left( q_{1}^{\ast }\right) .\Box \medskip $

\textbf{Proposition 2.4.} \textit{Let }$x,y\in \mathbb{O,}$ \textit{be two
real octonions. The following relations hold:}

\textit{1)} $iy=y^{\ast }i;$

\textit{2)} $\left( iy\right) x=i\left( y^{\ast }x^{\ast }\right) ^{\ast };$

\textit{3)} $x\left( iy\right) =i(x^{\ast }y);$

\textit{4)} $\left( iy\right) \left( ix\right) =-\left( yx^{\ast }\right)
^{\ast }.\medskip \medskip $ \ \ \ \ \ 

\textbf{Proof.}\newline
\qquad 1) We have $iy=e_{1}\left(
y_{0}+y_{1}e_{1}+y_{2}e_{2}+y_{3}e_{3}+y_{4}e_{4}+y_{5}e_{5}+y_{6}e_{6}+y_{7}e_{7}\right) = 
$\newline
$=$ $%
-y_{1}+y_{0}e_{1}-y_{3}e_{2}+y_{2}e_{3}-y_{5}e_{4}+y_{4}e_{5}+y_{7}e_{6}-y_{6}e_{7} 
$ and\newline
$y^{\ast }i=\left(
y_{0}+y_{1}e_{1}-y_{2}e_{2}-y_{3}e_{3}-y_{4}e_{4}-y_{5}e_{5}-y_{6}e_{6}-y_{7}e_{7}\right) e_{1}= 
$\newline
$%
=-y_{1}+y_{0}e_{1}-y_{3}e_{2}+y_{2}e_{3}-y_{5}e_{4}+y_{4}e_{5}+y_{7}e_{6}-y_{6}e_{7}. 
$

2) From the above, it results $\left( iy\right) x=\left( y^{\ast }i\right)
x. $ We have $\left( \left( iy\right) x\right) i=\left( \left( y^{\ast
}i\right) x\right) i.$ We apply relation $\left( 2.2\right) $ and we obtain $%
\left( \left( iy\right) x\right) i=y^{\ast }\left( ixi\right) =-$ $y^{\ast
}x^{\ast }.$ Therefore $\left( \left( \left( iy\right) x\right) i\right)
i=\left( -y^{\ast }x^{\ast }\right) i,$ then, from alternativity, we get $%
\left( iy\right) x=\left( y^{\ast }x^{\ast }\right) i=i\left( y^{\ast
}x^{\ast }\right) ^{\ast }.$

3) From relation $\left( 2.1\right) ,\,\ $the following relation holds 
\begin{equation*}
i(x\left( iy\right) )=\left( ixi\right) y=-x^{\ast }y.
\end{equation*}%
Using again alternativity, it results $i\left( i(x\left( iy\right) )\right)
=-i(x^{\ast }y),$ that means $x\left( iy\right) =i(x^{\ast }y).$

4) We apply relation $\left( 2.3\right) ~$and we have $\left( iy\right)
\left( ix\right) =\left( iy\right) \left( x^{\ast }i\right) =i\left(
yx^{\ast }\right) i=-\left( yx^{\ast }\right) ^{\ast }$.$\medskip ~\Box
\smallskip $

\textbf{Proposition 2.5. }\textit{For real octonion} $a=q_{1}+q_{2}k,$ 
\textit{with} $q_{1},q_{2}$ \textit{two real quaternions, we have the
following relations:}%
\begin{equation*}
\varepsilon \Lambda \left( a\right) \varepsilon =\Lambda \left( \widetilde{a}%
\right) ,\varepsilon \Delta \left( a\right) \varepsilon =\Delta \left( 
\widetilde{a}\right) 
\end{equation*}%
\textit{and}%
\begin{equation*}
\tau \Lambda \left( a\right) \tau =\Lambda \left( a_{+}\right) ,\tau \Delta
\left( a\right) \tau =\Delta \left( a_{+}\right) .
\end{equation*}

\bigskip

\textbf{Proof.} 1) Since $\varepsilon =diag(1,1,1,1,-1,-1,-1,-1)\in \mathcal{%
M}_{8}\left( \mathbb{R}\right) ,$ we have $\varepsilon =\left( 
\begin{array}{cc}
I_{4} & O_{4} \\ 
O_{4} & -I_{4}%
\end{array}%
\right) $, where $I_{4}\in \mathcal{M}_{4}\left( \mathbb{R}\right) $ is the
unit matrix and $O_{4}\in \mathcal{M}_{4}\left( \mathbb{R}\right) $ is the
zero matrix.

$\varepsilon \Lambda \left( a\right) \varepsilon =\left( 
\begin{array}{cc}
I_{4} & 0 \\ 
0 & -I_{4}%
\end{array}%
\right) \left( 
\begin{array}{cc}
\lambda \left( q_{1}\right) & -\rho \left( q_{2}\right) M_{1} \\ 
\lambda \left( q_{2}\right) M_{1} & \rho \left( q_{1}\right)%
\end{array}%
\right) \left( 
\begin{array}{cc}
I_{4} & 0 \\ 
0 & -I_{4}%
\end{array}%
\right) =$\newline
$=\left( 
\begin{array}{cc}
\lambda \left( q_{1}\right) & -\rho \left( q_{2}\right) M_{1} \\ 
-\lambda \left( q_{2}\right) M_{1} & -\rho \left( q_{1}\right)%
\end{array}%
\right) \left( 
\begin{array}{cc}
I_{4} & 0 \\ 
0 & -I_{4}%
\end{array}%
\right) =\left( 
\begin{array}{cc}
\lambda \left( q_{1}\right) & \rho \left( q_{2}\right) M_{1} \\ 
-\lambda \left( q_{2}\right) M_{1} & \rho \left( q_{1}\right)%
\end{array}%
\right) =$\newline
$=\Lambda \left( \widetilde{a}\right) .$

We have $\varepsilon \Delta \left( a\right) \varepsilon =\left( 
\begin{array}{cc}
I_{4} & 0 \\ 
0 & -I_{4}%
\end{array}%
\right) \left( 
\begin{array}{cc}
\rho \left( q_{1}\right) & -\lambda \left( \overline{q}_{2}\right) \\ 
\lambda \left( q_{2}\right) & \rho \left( \overline{q}_{1}\right)%
\end{array}%
\right) \left( 
\begin{array}{cc}
I_{4} & 0 \\ 
0 & -I_{4}%
\end{array}%
\right) =$\newline
$=\left( 
\begin{array}{cc}
\rho \left( q_{1}\right) & -\lambda \left( \overline{q}_{2}\right) \\ 
-\lambda \left( q_{2}\right) & -\rho \left( \overline{q}_{1}\right)%
\end{array}%
\right) \left( 
\begin{array}{cc}
I_{4} & 0 \\ 
0 & -I_{4}%
\end{array}%
\right) =\left( 
\begin{array}{cc}
\rho \left( q_{1}\right) & \lambda \left( \overline{q}_{2}\right) \\ 
-\lambda \left( q_{2}\right) & \rho \left( \overline{q}_{1}\right)%
\end{array}%
\right) =\Delta \left( \widetilde{a}\right) .$

$\bigskip $2) Since $\tau =diag(1,1,-1,-1,1,1,-1,-1),$ and $\sigma
=diag(1,1,-1,-1),~$we have\newline
$\tau \Lambda \left( a\right) \tau =\left( 
\begin{array}{cc}
\sigma & 0 \\ 
0 & \sigma%
\end{array}%
\right) \left( 
\begin{array}{cc}
\lambda \left( q_{1}\right) & -\rho \left( q_{2}\right) M_{1} \\ 
\lambda \left( q_{2}\right) M_{1} & \rho \left( q_{1}\right)%
\end{array}%
\right) \left( 
\begin{array}{cc}
\sigma & 0 \\ 
0 & \sigma%
\end{array}%
\right) =$\newline
$=\left( 
\begin{array}{cc}
\sigma \lambda \left( q_{1}\right) & -\sigma \rho \left( q_{2}\right) M_{1}
\\ 
\sigma \lambda \left( q_{2}\right) M_{1} & \sigma \rho \left( q_{1}\right)%
\end{array}%
\right) \left( 
\begin{array}{cc}
\sigma & 0 \\ 
0 & \sigma%
\end{array}%
\right) =$\newline
$=\left( 
\begin{array}{cc}
\sigma \lambda \left( q_{1}\right) \sigma & -\sigma \rho \left( q_{2}\right)
M_{1}\sigma \\ 
\sigma \lambda \left( q_{2}\right) M_{1}\sigma & \sigma \rho \left(
q_{1}\right) \sigma%
\end{array}%
\right) =\left( 
\begin{array}{cc}
\lambda \left( q_{1}^{\ast }\right) & -\rho \left( q_{2}^{\ast }\right) M_{1}
\\ 
\lambda \left( q_{2}^{\ast }\right) M_{1} & \rho \left( q_{1}^{\ast }\right)%
\end{array}%
\right) =\Lambda \left( a_{+}\right) ,~$\newline
since $M_{1}\sigma =\sigma M_{1}.$

$\tau \Delta \left( a\right) \tau =\left( 
\begin{array}{cc}
\sigma & 0 \\ 
0 & \sigma%
\end{array}%
\right) \left( 
\begin{array}{cc}
\rho \left( q_{1}\right) & -\lambda \left( \overline{q}_{2}\right) \\ 
\lambda \left( q_{2}\right) & \rho \left( \overline{q}_{1}\right)%
\end{array}%
\right) \left( 
\begin{array}{cc}
\sigma & 0 \\ 
0 & \sigma%
\end{array}%
\right) =$\newline
$=\left( 
\begin{array}{cc}
\sigma \rho \left( q_{1}\right) & -\sigma \lambda \left( \overline{q}%
_{2}\right) \\ 
\sigma \lambda \left( q_{2}\right) & \sigma \rho \left( \overline{q}%
_{1}\right)%
\end{array}%
\right) \left( 
\begin{array}{cc}
\sigma & 0 \\ 
0 & \sigma%
\end{array}%
\right) =\left( 
\begin{array}{cc}
\sigma \rho \left( q_{1}\right) \sigma & -\sigma \lambda \left( \overline{q}%
_{2}\right) \sigma \\ 
\sigma \lambda \left( q_{2}\right) \sigma & \sigma \rho \left( \overline{q}%
_{1}\right) \sigma%
\end{array}%
\right) =$\newline
$=\left( 
\begin{array}{cc}
\rho \left( q_{1}^{\ast }\right) & -\lambda \left( \overline{q}_{2}^{\ast
}\right) \\ 
\lambda \left( q_{2}^{\ast }\right) & \rho \left( \overline{q}_{1}^{\ast
}\right)%
\end{array}%
\right) =\Delta \left( a_{+}\right) $.$~\Box \medskip $

\textbf{Proposition 2.6.} \textit{Let} $A,X\in \mathbb{O}_{\mathbb{C}%
},A=x+iy,X=v+iw,x,y,v,w\in \mathbb{O},$ \textit{then:}

\begin{equation*}
\overrightarrow{AX}=\Phi \left( A\right) \overrightarrow{X}.
\end{equation*}

\textbf{Proof.}

We have\newline
\begin{equation*}
\Phi \left( A\right) \overrightarrow{X}=\bigskip \left( 
\begin{array}{cc}
\Lambda \left( x\right) & -\Lambda \left( y\right) \\ 
\Lambda \left( y^{\ast }\right) & \Lambda \left( x^{\ast }\right)%
\end{array}%
\right) \left( 
\begin{array}{c}
\overrightarrow{v} \\ 
\overrightarrow{w}%
\end{array}%
\right) =
\end{equation*}%
\newline
\begin{equation*}
=\left( 
\begin{array}{c}
\Lambda \left( x\right) \overrightarrow{v}-\Lambda \left( y\right) 
\overrightarrow{w} \\ 
\Lambda \left( y^{\ast }\right) \overrightarrow{v}+\Lambda \left( x^{\ast
}\right) \overrightarrow{w}%
\end{array}%
\right) .
\end{equation*}

We have $AX=\left( x+iy\right) \left( v+iw\right) =xv+x\left( iw\right)
+\left( iy\right) v+\left( iy\right) \left( iw\right) =$\newline
$=xv-\left( yw^{\ast }\right) ^{\ast }+i(x^{\ast }w+\left( y^{\ast }v^{\ast
}\right) ^{\ast })$. From Proposition 2.4, i) and relations 2.22, 2.23 and
2.24, it follows that 
\begin{equation*}
\overrightarrow{\left( yw^{\ast }\right) ^{\ast }}=-\overrightarrow{i\left(
yw^{\ast }\right) i}=-\Lambda \left( i\right) \Delta \left( i\right) 
\overrightarrow{yw^{\ast }}=
\end{equation*}%
\begin{equation*}
=-\Lambda \left( i\right) \Delta \left( i\right) \Lambda \left( y\right) 
\overrightarrow{w^{\ast }}=\Lambda \left( i\right) \Delta \left( i\right)
\Lambda \left( y\right) \overrightarrow{iwi}=
\end{equation*}%
\begin{equation*}
=\Lambda \left( i\right) \Delta \left( i\right) \Lambda \left( y\right)
\Lambda \left( i\right) \Delta \left( i\right) \overrightarrow{w}=
\end{equation*}%
\begin{equation*}
\Lambda \left( i^{2}\right) \Delta \left( i^{2}\right) \Lambda \left(
y\right) \overrightarrow{w}=\Lambda \left( y\right) \overrightarrow{w},
\end{equation*}%
since $\Lambda \left( i^{2}\right) =\Delta \left( i^{2}\right) =\Lambda
\left( -1\right) =\Delta \left( -1\right) =-I_{8}.$ Therefore, we have 
\begin{equation*}
\overrightarrow{AX}=\left( 
\begin{array}{c}
\overrightarrow{xv}-\overrightarrow{\left( yw^{\ast }\right) ^{\ast }} \\ 
\overrightarrow{x^{\ast }w}+\overrightarrow{\left( y^{\ast }v^{\ast }\right)
^{\ast }}%
\end{array}%
\right) =\left( 
\begin{array}{c}
\overrightarrow{xv}-\Lambda \left( y\right) \overrightarrow{w} \\ 
\overrightarrow{x^{\ast }w}+\Lambda \left( y^{\ast }\right) \overrightarrow{v%
}%
\end{array}%
\right) =
\end{equation*}

\begin{equation*}
=\left( 
\begin{array}{c}
\Lambda \left( x\right) \overrightarrow{v}-\Lambda \left( y\right) 
\overrightarrow{w} \\ 
\Lambda \left( y^{\ast }\right) \overrightarrow{v}+\Lambda \left( x^{\ast
}\right) \overrightarrow{w}%
\end{array}%
\right) =
\end{equation*}

\begin{equation*}
=\left( 
\begin{array}{cc}
\Lambda \left( x\right) & -\Lambda \left( y\right) \\ 
\Lambda \left( y^{\ast }\right) & \Lambda \left( x^{\ast }\right)%
\end{array}%
\right) \left( 
\begin{array}{c}
\overrightarrow{v} \\ 
\overrightarrow{w}%
\end{array}%
\right) =\Phi \left( A\right) \overrightarrow{X}.\Box \medskip
\end{equation*}

Let $M,N$ be the matrices 
\begin{equation*}
N=\left( 1,e_{1},e_{2},e_{3},e_{4},e_{5},e_{6},e_{7}\right) ^{t},
\end{equation*}%
\begin{equation*}
M=\left( 1,-e_{1},-e_{2},-e_{3},-e_{4},-e_{5},-e_{6},-e_{7}\right) ^{t}.
\end{equation*}%
We remark that $\medskip N^{t}M=8.\medskip $

\textbf{Proposition 2.7.} \textit{If \ }$A\in \mathbb{O},$ $%
a=a_{0}+a_{1}e_{1}+a_{2}e_{2}+a_{3}e_{3}+a_{4}e_{4}+a_{5}e_{5}+a_{6}e_{6}+a_{7}e_{7} 
$ \textit{we have:}

\textit{i)} \ $\Lambda \left( a\right) M=Ma.$

\textit{ii)} $\theta M=Me_{1}.$

\textit{iii)} $\Lambda \left( a\right) N=\overline{a}N.\medskip $

\textbf{Proof.} i) We have $\Lambda \left( a\right) M=$\newline
\begin{equation*}
=\left( 
\begin{array}{cccccccc}
a_{0} & -a_{1} & -a_{2} & -a_{3} & -a_{4} & -a_{5} & -a_{6} & -a_{7} \\ 
a_{1} & a_{0} & -a_{3} & a_{2} & -a_{5} & a_{4} & a_{7} & -a_{6} \\ 
a_{2} & a_{3} & a_{0} & -a_{1} & -a_{6} & -a_{7} & a_{4} & a_{5} \\ 
a_{3} & -a_{2} & a_{1} & a_{0} & -a_{7} & a_{6} & -a_{5} & a_{4} \\ 
a_{4} & a_{5} & a_{6} & a_{7} & a_{0} & -a_{1} & -a_{2} & -a_{3} \\ 
a_{5} & -a_{4} & a_{7} & -a_{6} & a_{1} & a_{0} & a_{3} & -a_{2} \\ 
a_{6} & -a_{7} & -a_{4} & a_{5} & a_{2} & -a_{3} & a_{0} & a_{1} \\ 
a_{7} & a_{6} & -a_{5} & -a_{4} & a_{3} & a_{2} & -a_{1} & a_{0}%
\end{array}%
\right) \left( 
\begin{array}{c}
1 \\ 
-e_{1} \\ 
-e_{2} \\ 
-e_{3} \\ 
-e_{4} \\ 
-e_{5} \\ 
-e_{6} \\ 
-e_{7}%
\end{array}%
\right) =
\end{equation*}

\begin{equation*}
=\left( 
\begin{array}{c}
a_{0}+a_{1}e_{1}+a_{2}e_{2}+a_{3}e_{3}+a_{4}e_{4}+a_{5}e_{5}+a_{6}e_{6}+a_{7}e_{7}
\\ 
a_{1}-a_{0}e_{1}+a_{3}e_{2}-a_{2}e_{3}+a_{5}e_{4}-a_{4}e_{5}-a_{7}e_{6}+a_{6}e_{7}
\\ 
a_{2}-a_{3}e_{1}-a_{0}e_{2}+a_{1}e_{3}+a_{6}e_{4}+a_{7}e_{5}-a_{4}e_{6}-a_{5}e_{7}
\\ 
a_{3}+a_{2}e_{1}-a_{1}e_{2}-a_{0}e_{3}+a_{7}e_{4}-a_{6}e_{5}+a_{5}e_{6}-a_{4}e_{7}
\\ 
a_{4}-a_{5}e_{1}-a_{6}e_{2}-a_{7}e_{3}-a_{0}e_{4}+a_{1}e_{5}+a_{2}e_{6}+a_{3}e_{7}
\\ 
a_{5}+a_{4}e_{1}-a_{7}e_{2}+a_{6}e_{3}-a_{1}e_{4}-a_{0}e_{5}-a_{3}e_{6}+a_{2}e_{7}
\\ 
a_{6}+a_{7}e_{1}+a_{4}e_{2}-a_{5}e_{3}-a_{2}e_{4}+a_{3}e_{5}-a_{0}e_{6}-a_{1}e_{7}
\\ 
a_{7}-a_{6}e_{1}+a_{5}e_{2}+a_{4}e_{3}-a_{3}e_{4}-a_{2}e_{5}+a_{1}e_{6}-a_{0}e_{7}%
\end{array}%
\right) =
\end{equation*}

\begin{equation*}
=\left( 
\begin{array}{c}
a_{0}+a_{1}e_{1}+a_{2}e_{2}+a_{3}e_{3}+a_{4}e_{4}+a_{5}e_{5}+a_{6}e_{6}+a_{7}e_{7}
\\ 
-e_{1}\left(
a_{0}+a_{1}e_{1}+a_{2}e_{2}+a_{3}e_{3}+a_{4}e_{4}+a_{5}e_{5}+a_{6}e_{6}+a_{7}e_{7}\right)
\\ 
-e_{2}(a_{0}+a_{1}e_{1}+a_{2}e_{2}+a_{3}e_{3}+a_{4}e_{4}+a_{5}e_{5}+a_{6}e_{6}+a_{7}e_{7})
\\ 
-e_{3}(a_{0}+a_{1}e_{1}+a_{2}e_{2}+a_{3}e_{3}+a_{4}e_{4}+a_{5}e_{5}+a_{6}e_{6}+a_{7}e_{7})
\\ 
-e_{4}(a_{0}+a_{1}e_{1}+a_{2}e_{2}+a_{3}e_{3}+a_{4}e_{4}+a_{5}e_{5}+a_{6}e_{6}+a_{7}e_{7})
\\ 
-e_{5}(a_{0}+a_{1}e_{1}+a_{2}e_{2}+a_{3}e_{3}+a_{4}e_{4}+a_{5}e_{5}+a_{6}e_{6}+a_{7}e_{7})
\\ 
-e_{6}(a_{0}+a_{1}e_{1}+a_{2}e_{2}+a_{3}e_{3}+a_{4}e_{4}+a_{5}e_{5}+a_{6}e_{6}+a_{7}e_{7})
\\ 
-e_{7}(a_{0}+a_{1}e_{1}+a_{2}e_{2}+a_{3}e_{3}+a_{4}e_{4}+a_{5}e_{5}+a_{6}e_{6}+a_{7}e_{7})%
\end{array}%
\right) =
\end{equation*}

$=Ma.$

ii) 
\begin{equation*}
\theta M=\left( 
\begin{array}{cccccccc}
0 & -1 & 0 & 0 & 0 & 0 & 0 & 0 \\ 
1 & 0 & 0 & 0 & 0 & 0 & 0 & 0 \\ 
0 & 0 & 0 & -1 & 0 & 0 & 0 & 0 \\ 
0 & 0 & 1 & 0 & 0 & 0 & 0 & 0 \\ 
0 & 0 & 0 & 0 & 0 & -1 & 0 & 0 \\ 
0 & 0 & 0 & 0 & 1 & 0 & 0 & 0 \\ 
0 & 0 & 0 & 0 & 0 & 0 & 0 & 1 \\ 
0 & 0 & 0 & 0 & 0 & 0 & -1 & 0%
\end{array}%
\right) \left( 
\begin{array}{c}
1 \\ 
-e_{1} \\ 
-e_{2} \\ 
-e_{3} \\ 
-e_{4} \\ 
-e_{5} \\ 
-e_{6} \\ 
-e_{7}%
\end{array}%
\right) =
\end{equation*}

\begin{equation*}
=\left( 
\begin{array}{c}
e_{1} \\ 
1 \\ 
e_{3} \\ 
-e_{2} \\ 
e_{5} \\ 
-e_{4} \\ 
-e_{7} \\ 
e_{6}%
\end{array}%
\right) =Me_{1}.
\end{equation*}

iii) 
\begin{equation*}
\Lambda \left( a\right) N=\left( 
\begin{array}{cccccccc}
a_{0} & -a_{1} & -a_{2} & -a_{3} & -a_{4} & -a_{5} & -a_{6} & -a_{7} \\ 
a_{1} & a_{0} & -a_{3} & a_{2} & -a_{5} & a_{4} & a_{7} & -a_{6} \\ 
a_{2} & a_{3} & a_{0} & -a_{1} & -a_{6} & -a_{7} & a_{4} & a_{5} \\ 
a_{3} & -a_{2} & a_{1} & a_{0} & -a_{7} & a_{6} & -a_{5} & a_{4} \\ 
a_{4} & a_{5} & a_{6} & a_{7} & a_{0} & -a_{1} & -a_{2} & -a_{3} \\ 
a_{5} & -a_{4} & a_{7} & -a_{6} & a_{1} & a_{0} & a_{3} & -a_{2} \\ 
a_{6} & -a_{7} & -a_{4} & a_{5} & a_{2} & -a_{3} & a_{0} & a_{1} \\ 
a_{7} & a_{6} & -a_{5} & -a_{4} & a_{3} & a_{2} & -a_{1} & a_{0}%
\end{array}%
\right) \left( 
\begin{array}{c}
1 \\ 
e_{1} \\ 
e_{2} \\ 
e_{3} \\ 
e_{4} \\ 
e_{5} \\ 
e_{6} \\ 
e_{7}%
\end{array}%
\right) =\newline
\end{equation*}%
\newline
\begin{equation*}
=\left( 
\begin{array}{c}
a_{0}-a_{1}e_{1}-a_{2}e_{2}-a_{3}e_{3}-a_{4}e_{4}-a_{5}e_{5}-a_{6}e_{6}-a_{7}e_{7}
\\ 
a_{1}+a_{0}e_{1}-a_{3}e_{2}+a_{2}e_{3}-a_{5}e_{4}+a_{4}e_{5}+a_{7}e_{6}-a_{6}e_{7}
\\ 
a_{2}+a_{3}e_{1}+a_{0}e_{2}-a_{1}e_{3}-a_{6}e_{4}-a_{7}e_{5}+a_{4}e_{6}+a_{5}e_{7}
\\ 
a_{3}-a_{2}e_{1}+a_{1}e_{2}+a_{0}e_{3}-a_{7}e_{4}+a_{6}e_{5}-a_{5}e_{6}+a_{4}e_{7}
\\ 
a_{4}+a_{5}e_{1}+a_{6}e_{2}+a_{7}e_{3}+a_{0}e_{4}-a_{1}e_{5}-a_{2}e_{6}-a_{3}e_{7}
\\ 
a_{5}-a_{4}e_{1}+a_{7}e_{2}-a_{6}e_{3}+a_{1}e_{4}+a_{0}e_{5}+a_{3}e_{6}-a_{2}e_{7}
\\ 
a_{6}-a_{7}e_{1}-a_{4}e_{2}+a_{5}e_{3}+a_{2}e_{4}-a_{3}e_{5}+a_{0}e_{6}+a_{1}e_{7}
\\ 
a_{7}+a_{6}e_{1}-a_{5}e_{2}-a_{4}e_{3}+a_{3}e_{4}+a_{2}e_{5}-a_{1}e_{6}+a_{0}e_{7}%
\end{array}%
\right) =\newline
\end{equation*}%
\begin{equation*}
=\left( 
\begin{array}{c}
a_{0}-a_{1}e_{1}-a_{2}e_{2}-a_{3}e_{3}-a_{4}e_{4}-a_{5}e_{5}-a_{6}e_{6}-a_{7}e_{7}
\\ 
\left(
a_{0}-a_{1}e_{1}-a_{2}e_{2}-a_{3}e_{3}-a_{4}e_{4}-a_{5}e_{5}-a_{6}e_{6}-a_{7}e_{7}\right) e_{1}
\\ 
\left(
a_{0}-a_{1}e_{1}-a_{2}e_{2}-a_{3}e_{3}-a_{4}e_{4}-a_{5}e_{5}-a_{6}e_{6}-a_{7}e_{7}\right) e_{2}
\\ 
\left(
a_{0}-a_{1}e_{1}-a_{2}e_{2}-a_{3}e_{3}-a_{4}e_{4}-a_{5}e_{5}-a_{6}e_{6}-a_{7}e_{7}\right) e_{3}
\\ 
\left(
a_{0}-a_{1}e_{1}-a_{2}e_{2}-a_{3}e_{3}-a_{4}e_{4}-a_{5}e_{5}-a_{6}e_{6}-a_{7}e_{7}\right) e_{4}
\\ 
\left(
a_{0}-a_{1}e_{1}-a_{2}e_{2}-a_{3}e_{3}-a_{4}e_{4}-a_{5}e_{5}-a_{6}e_{6}-a_{7}e_{7}\right) e_{5}
\\ 
\left(
a_{0}-a_{1}e_{1}-a_{2}e_{2}-a_{3}e_{3}-a_{4}e_{4}-a_{5}e_{5}-a_{6}e_{6}-a_{7}e_{7}\right) e_{6}
\\ 
\left(
a_{0}-a_{1}e_{1}-a_{2}e_{2}-a_{3}e_{3}-a_{4}e_{4}-a_{5}e_{5}-a_{6}e_{6}-a_{7}e_{7}\right) e_{7}%
\end{array}%
\right) \newline
=
\end{equation*}%
$=\overline{a}N.\Box \medskip $

\textbf{Proposition 2.8.} \textit{Let} $A=x+iy$ \textit{be a complex
octonion with} $x,y$ \textit{two real octonions and }$a=q_{1}+q_{2}k$ 
\textit{be a real octonion, with} $q_{1},q_{2}$ \textit{two real
quaternions. The following relations are true.}

\textit{i)} $T\Lambda \left( a\right) T=\Delta \left( a^{+}\right) ,$\textit{%
where} $T=\left( 
\begin{array}{cc}
M_{1} & O_{4} \\ 
O_{4} & I_{4}%
\end{array}%
\right) \in \mathcal{M}_{8}\left( \mathbb{R}\right) $ \textit{and} $O_{4}\in 
\mathcal{M}_{8}\left( \mathbb{R}\right) $ \textit{is zero matrix.}

\textit{ii)} $S\Phi \left( A\right) S=\Psi \left( A^{+}\right) ,$ \textit{%
where} $A^{+}=x^{+}+iy^{+},S=\left( 
\begin{array}{cc}
T & O_{8} \\ 
O_{8} & T%
\end{array}%
\right) \in \mathcal{M}_{16}\left( \mathbb{R}\right) ~$\textit{and} $%
O_{8}\in \mathcal{M}_{8}\left( \mathbb{R}\right) ~$\textit{is zero matrix.}%
\medskip

\textbf{Proof.} i) From (Tian, 2000), relations 1.17 and 1.18, we know that $%
\rho \left( q\right) =M_{1}\lambda ^{t}\left( q\right) M_{1}=M_{1}\lambda
\left( \overline{q}\right) M_{1},$ with $q$ a real quaternion$.$ It results 
\newline
\begin{equation*}
\left( 
\begin{array}{cc}
M_{1} & O_{4} \\ 
O_{4} & I_{4}%
\end{array}%
\right) \Lambda \left( A\right) \left( 
\begin{array}{cc}
M_{1} & O_{4} \\ 
O_{4} & I_{4}%
\end{array}%
\right) =\newline
\end{equation*}%
\begin{equation*}
=\left( 
\begin{array}{cc}
M_{1} & O_{4} \\ 
O_{4} & I_{4}%
\end{array}%
\right) \left( 
\begin{array}{cc}
\lambda \left( q_{1}\right) & -\rho \left( q_{2}\right) M_{1} \\ 
\lambda \left( q_{2}\right) M_{1} & \rho \left( q_{1}\right)%
\end{array}%
\right) \left( 
\begin{array}{cc}
M_{1} & O_{4} \\ 
O_{4} & I_{4}%
\end{array}%
\right) =\newline
\end{equation*}%
\begin{equation*}
=\left( 
\begin{array}{cc}
M_{1}\lambda \left( q_{1}\right) & -M_{1}\rho \left( q_{2}\right) M_{1} \\ 
\lambda \left( q_{2}\right) M_{1} & \rho \left( q_{1}\right)%
\end{array}%
\right) \left( 
\begin{array}{cc}
M_{1} & O_{4} \\ 
O_{4} & I_{4}%
\end{array}%
\right) =\newline
\end{equation*}%
\begin{equation*}
=\left( 
\begin{array}{cc}
M_{1}\lambda \left( q_{1}\right) M_{1} & -M_{1}\rho \left( q_{2}\right) M_{1}
\\ 
\lambda \left( q_{2}\right) & \rho \left( q_{1}\right)%
\end{array}%
\right) =\left( 
\begin{array}{cc}
\rho \left( \overline{q}_{1}\right) & -M_{1}\rho \left( q_{2}\right) M_{1}
\\ 
\lambda \left( q_{2}\right) & \rho \left( q_{1}\right)%
\end{array}%
\right) =
\end{equation*}

\begin{equation*}
=\left( 
\begin{array}{cc}
\rho \left( \overline{q}_{1}\right) & -\lambda \left( \overline{q}_{2}\right)
\\ 
\lambda \left( q_{2}\right) & \rho \left( q_{1}\right)%
\end{array}%
\right) =\Delta \left( a^{+}\right) .
\end{equation*}

\bigskip ii) We have 
\begin{equation*}
\left( 
\begin{array}{cc}
T & O_{8} \\ 
O_{8} & T%
\end{array}%
\right) \Phi \left( A\right) \left( 
\begin{array}{cc}
T & O_{8} \\ 
O_{8} & T%
\end{array}%
\right) =\newline
\end{equation*}%
\begin{equation*}
=\left( 
\begin{array}{cc}
T & O_{8} \\ 
O_{8} & T%
\end{array}%
\right) \bigskip \left( 
\begin{array}{cc}
\Lambda \left( x\right) & -\Lambda \left( y\right) \\ 
\Lambda \left( y^{\ast }\right) & \Lambda \left( x^{\ast }\right)%
\end{array}%
\right) \left( 
\begin{array}{cc}
T & O_{8} \\ 
O_{8} & T%
\end{array}%
\right) =\newline
\end{equation*}%
\begin{equation*}
=\left( 
\begin{array}{cc}
T\Lambda \left( x\right) & -T\Lambda \left( y\right) \\ 
T\Lambda \left( y^{\ast }\right) & T\Lambda \left( x^{\ast }\right)%
\end{array}%
\right) \left( 
\begin{array}{cc}
T & O_{8} \\ 
O_{8} & T%
\end{array}%
\right) =\left( 
\begin{array}{cc}
T\Lambda \left( x\right) T & -T\Lambda \left( y\right) T \\ 
T\Lambda \left( y^{\ast }\right) T & T\Lambda \left( x^{\ast }\right) T%
\end{array}%
\right) =\newline
\end{equation*}%
\begin{equation*}
=\left( 
\begin{array}{cc}
\Delta \left( x^{+}\right) & -\Delta \left( y^{+}\right) \\ 
\Delta \left( y^{\ast +}\right) & \Delta \left( x^{\ast +}\right)%
\end{array}%
\right) =\Psi \left( A^{+}\right) .\Box \medskip
\end{equation*}

\textbf{Proposition 2.9.} \textit{With the above notations, we have}

\begin{equation*}
\Psi \left( A^{+}\right) \overrightarrow{X}=S\overrightarrow{AX^{+}},
\end{equation*}%
\textit{where} $X^{+}=v^{+}+iw^{+}.$

\medskip

\textbf{Proof.} Since $T\overrightarrow{v}=\overrightarrow{v^{+}},$ it
results \newline
$\Psi \left( A^{+}\right) \overrightarrow{X}=$\bigskip $S\Phi \left(
A\right) S\overrightarrow{X}=S\Phi \left( A\right) \left( 
\begin{array}{cc}
T & O_{8} \\ 
O_{8} & T%
\end{array}%
\right) \left( 
\begin{array}{c}
\overrightarrow{v} \\ 
\overrightarrow{w}%
\end{array}%
\right) =$\newline
$=S\Phi \left( A\right) \left( 
\begin{array}{c}
T\overrightarrow{v} \\ 
T\overrightarrow{w}%
\end{array}%
\right) =S\Phi \left( A\right) \left( 
\begin{array}{c}
\overrightarrow{v^{+}} \\ 
\overrightarrow{w^{+}}%
\end{array}%
\right) =S\Phi \left( A\right) \overrightarrow{X^{+}}=S\overrightarrow{AX^{+}%
}.\Box $

\begin{equation*}
\end{equation*}%
\begin{equation*}
\end{equation*}

\textbf{3. A set of invertible elements in split quaternion and octonion
algebras}%
\begin{equation*}
\end{equation*}

In a split quaternion algebra and in a split octonion algebra there are
nonzero elements such that their norms are zero. In such algebras, it is
very good to know sets of invertible elements, that means nonzero elements
with their norms nonzero. In the following, we give a method to find such a
sets.

Let $n$ be an arbitrary positive integer and let $a,b,c,x_{0},x_{1},x_{2}~$%
be arbitrary integers. We consider the following difference equation of
degree three%
\begin{equation}
X_{n}=aX_{n-1}+bX_{n-2}+cX_{n-3},X_{0}=x_{0},X_{1}=x_{1},X_{2}=x_{2}. 
\tag{3.1}
\end{equation}

We consider the following degree three equation%
\begin{equation}
x^{3}-ax^{2}-bx-c=0.  \tag{3.2}
\end{equation}%
We consider that this equation has three real solutions $\sigma _{1}>\sigma
_{2}>\sigma _{3},$ with $\sigma _{1}>1.~$For this case, we have the
following Binet's formula:%
\begin{equation}
X_{n}=A\sigma _{1}^{n}+B\sigma _{2}^{n}+C\sigma _{3}^{n},  \tag{3.3}
\end{equation}%
where $A,B,C$ are solutions of the following linear system:%
\begin{equation*}
\left\{ 
\begin{array}{c}
A+B+C=x_{0} \\ 
A\sigma _{1}+B\sigma _{2}+C\sigma _{3}=x_{1} \\ 
A\sigma _{1}^{2}+B\sigma _{2}^{2}+C\sigma _{3}^{2}=x_{2}%
\end{array}%
\right. .
\end{equation*}%
Since we obtain a Vandermonde determinant, the system has a unique solution.

We consider the real generalized quaternion algebra $\mathbb{H}\left( \beta
_{1},\beta _{2}\right) $ and we define the quaternions 
\begin{equation*}
W_{n}=X_{n}+X_{n+1}e_{2}+X_{n+2}e_{3}+X_{n+3}e_{4},
\end{equation*}%
where $X_{n}$ is the $n$th number given by the relation $\left( 3.3\right) .$

From the above, we can compute the following limit 
\begin{equation*}
\underset{n\rightarrow \infty }{\lim }\boldsymbol{n}\left( W_{n}\right) =%
\underset{n\rightarrow \infty }{\lim }(X_{n}^{2}+\beta _{1}X_{n+1}^{2}+\beta
_{2}X_{n+2}^{2}+\beta _{1}\beta _{2}X_{n+3}^{2})=
\end{equation*}%
\begin{eqnarray*}
&=&\underset{n\rightarrow \infty }{\lim }((A\sigma _{1}^{n}+B\sigma
_{2}^{n}+C\sigma _{3}^{n})^{2}+\beta _{1}\left( A\sigma _{1}^{n+1}+B\sigma
_{2}^{n+1}+C\sigma _{3}^{n+1}\right) ^{2}+ \\
&&\beta _{2}\left( A\sigma _{1}^{n+2}+B\sigma _{2}^{n+2}+C\sigma
_{3}^{n+2}\right) ^{2}+\beta _{1}\beta _{2}\left( A\sigma _{1}^{n+3}+B\sigma
_{2}^{n+3}+C\sigma _{3}^{n+3}\right) ^{2}).
\end{eqnarray*}%
Let $f\left( \beta _{1},\beta _{2}\right) =A^{2}(1+\beta _{1}\sigma
_{1}^{2}+\beta _{2}\sigma _{1}^{4}+\beta _{1}\beta _{2}\sigma _{1}^{6}).$ If
\ $f\left( \beta _{1},\beta _{2}\right) \neq 0,$ it results that 
\begin{equation*}
\underset{n\rightarrow \infty }{\lim }\boldsymbol{n}\left( W_{n}\right)
=signf\left( \beta _{1},\beta _{2}\right) \cdot \infty .
\end{equation*}

Therefore, for \ all $\beta _{1},\beta _{2}\in \mathbb{R}$ with $f(\beta
_{1},\beta _{2})\neq 0,$ in the algebra $\mathbb{H}\left( \beta _{1},\beta
_{2}\right) $ there is a natural number $n_{0}$ such that $\boldsymbol{n}%
\left( W_{n}\right) \neq 0.$ From here, we have that $W_{n}$ is an
invertible element for all $n\geq n_{0}.$

Now, we consider the real octonion algebra $\mathbb{O}(\alpha ,\beta ,\gamma
).$ We define the octonions 
\begin{equation*}
Z_{n}=X_{n}+X_{n+1}e_{2}+X_{n+2}e_{3}+X_{n+3}e_{4}+X_{n+4}e_{5}+X_{n+5}e_{6}+X_{n+6}e_{7},
\end{equation*}%
where $X_{n}$ is the $n$th number given by the relation $\left( 3.3\right) .$

From the above, we can compute the following limit

\begin{equation*}
\underset{n\rightarrow \infty }{\lim }\boldsymbol{n}\left( Z_{n}\right) =
\end{equation*}%
\begin{equation*}
=\underset{n\rightarrow \infty }{\lim }(X_{n}^{2}\text{+}\alpha X_{n+1}^{2}%
\text{+}\beta X_{n+2}^{2}\text{+}\alpha \beta X_{n+3}^{2}\text{+}\gamma
X_{n+4}^{2}\text{+}\alpha \gamma X_{n+5}^{2}\text{+}\beta \gamma X_{n+6}^{2}%
\text{+}\alpha \beta \gamma X_{n+7}^{2})=
\end{equation*}%
\begin{eqnarray*}
&=&\underset{n\rightarrow \infty }{\lim }((A\sigma _{1}^{n}+B\sigma
_{2}^{n}+C\sigma _{3}^{n})^{2}+\alpha \left( A\sigma _{1}^{n+1}+B\sigma
_{2}^{n+1}+C\sigma _{3}^{n+1}\right) ^{2}+ \\
+ &&\beta \left( A\sigma _{1}^{n+2}+B\sigma _{2}^{n+2}+C\sigma
_{3}^{n+2}\right) ^{2}+\alpha \beta \left( A\sigma _{1}^{n+3}+B\sigma
_{2}^{n+3}+C\sigma _{3}^{n+3}\right) ^{2}+
\end{eqnarray*}%
\begin{equation*}
+\gamma \left( A\sigma _{1}^{n+4}+B\sigma _{2}^{n+4}+C\sigma
_{3}^{n+4}\right) ^{2}+\alpha \gamma \left( A\sigma _{1}^{n+5}+B\sigma
_{2}^{n+5}+C\sigma _{3}^{n+5}\right) ^{2}+
\end{equation*}%
\begin{equation*}
+\beta \gamma \left( A\sigma _{1}^{n+6}+B\sigma _{2}^{n+6}+C\sigma
_{3}^{n+6}\right) ^{2}+\alpha \beta \gamma \left( A\sigma _{1}^{n+7}+B\sigma
_{2}^{n+7}+C\sigma _{3}^{n+7}\right) ^{2}.
\end{equation*}

Let $g\left( \alpha ,\beta ,\gamma \right) $=$A^{2}(1$+$\alpha \sigma
_{1}^{2}$+$\beta \sigma _{1}^{4}$+$\alpha \beta \sigma _{1}^{6}$+$\gamma
\sigma _{1}^{8}$+$\alpha \gamma \sigma _{1}^{10}$+$\beta \gamma \sigma
_{1}^{12}$+$\alpha \beta \gamma \sigma _{1}^{14}).$

If \ $g\left( \alpha ,\beta ,\gamma \right) \neq 0,$ it results that 
\begin{equation*}
\underset{n\rightarrow \infty }{\lim }\boldsymbol{n}\left( Z_{n}\right)
=signg\left( \alpha ,\beta ,\gamma \right) \cdot \infty .
\end{equation*}

Therefore, for \ all $\alpha ,\beta ,\gamma \in \mathbb{R}$ with $g\left(
\alpha ,\beta ,\gamma \right) \neq 0,$ in the algebra $\mathbb{O}(\alpha
,\beta ,\gamma )$ there is a natural number $n_{0}$ such that $\boldsymbol{n}%
\left( Z_{n}\right) \neq 0.$ From here, we have that $z_{n}$ is an
invertible element for all $n\geq n_{0}.$

Since algebras $\mathbb{H}\left( \beta _{1},\beta _{2}\right) $ and $\mathbb{%
O}(\alpha ,\beta ,\gamma )~$are not always division algebras, finding
examples of invertible elements in such algebras can be a difficult problem.
The above elements, $W_{n}$ and $Z_{n},$ provide us \ an infinite set of
invertible elements in a split quaternion algebra and in a split octonion
algebra.%
\begin{equation*}
\end{equation*}

\textbf{Conclusions.} In this chapter, we gave some properties of the real
matrix representations for complex octonions and we provided sets of
invertible elements in a split quaternion algebra and in a split octonion
algebra. Due their applications, the study of these representations and the
study of these elements can give us other properties and applications.

\bigskip

\begin{equation*}
\end{equation*}%
\textbf{References}%
\begin{equation*}
\end{equation*}

Alamouti, S.M.: A simple transmit diversity technique for wireless
communications,\ IEEE J. Sel. Areas Commun., \textbf{16}(8), 1451-- 1458
(1998).

Alfsmann, D., G\"{o}ckler, H.G., Sangwine, S.J., Ell, T.,A.: Hypercomplex
Algebras in Digital Signal Processing: Benefits and Drawbacks, 15th European
Signal Processing Conference (EUSIPCO 2007), Poznan, Poland, 1322-1326
(2007).

Baez, J.C.: The Octonions, B. Am. Math. Soc., \textbf{39(}2\textbf{)},
145-205 (2002).

Chanyal, B.~C.: Octonion massive electrodynamics, General Relativity and
Gravitation, \textbf{46}, article ID: 1646, (2014).

Chanyal, B.~C., ~Bisht, P.~S., ~Negi, O.~P.~S.: Generalized Split-Octonion
Electrodynamics, International Journal of Theoretical Physics, \textbf{50}%
(6), 1919-1926 (2011).

Chen, J., Tu, A.: Fabric Image Edge Detection Based on Octonion and Echo
State Networks, Applied Mechanics and Matherials, \textbf{520}, 714-718 \
(2014).

Flaut, C.: BCK-algebras arising from block codes, J. Intell. Fuzzy Syst., 
\textbf{28}(4), 1829-1833 (2015).

Flaut, C.: Divison algebras with dimension 2\symbol{94}t, t in N, Analele
Stiintifice ale Universitatii Ovidius Constanta, Seria Matematica, \textbf{13%
}(2), 31-38 (2006).

Flaut, C., Savin, D.: Some examples of division symbol algebras of degree 3
and 5, Carpathian J. Math, \textbf{31}(2), 197-204 (2015).

Flaut, C., Savin, D.: About quaternion algebras and symbol algebras,
Bulletin of the Transilvania University of Brasov, 7(2)(56), 59-64 (2014).

Flaut, C., Savin, D., Iorgulecu, G.: Some properties of Fibonacci and Lucas
symbol elements, Journal of Mathematical Sciences: Advances and
Applications, 20, 37-43, (2013).

Flaut, C., Shpakivskyi, V.,: An Efficient Method for Solving Equations in
Generalized Quaternion and Octonion Algebras, Adv. Appl. Clifford Algebras,
25(2), 337-350 (2015).

Flaut, C., Shpakivskyi, V.,: Holomorphic functions in generalized Cayley-
Dickson algebras, Adv. Appl. Clifford Algebras, 25(1), 95-112 (2015).

Flaut, C., Shpakivskyi, V.: Real matrix representations for the complex
quaternions, Adv. Appl. Clifford Algebras, \textbf{23}(3), 657-671 (2013).

Flaut, C., Shpakivskyi, V.,: On Generalized Fibonacci Quaternions and
Fibonacci-Narayana Quaternions, Adv. Appl. Clifford Algebras, 23(3), 673-688
(2013).

Flaut, C., Shpakivskyi, V.,: Some Remarks About Fibonacci Elements in An
Arbitrary Algebra, Bulletin De La Soci\'{e}t\'{e} Des Sciences Et Des
Lettres De \L \'{o}dz, LXV(3), 63-73 (2015).

Flaut, C., \c{S}tefanescu, M.: Some equations over generalized quaternion
and octonion division algebras, Bull. Math. Soc. Sci. Math. Roumanie, 
\textbf{52}(4)(100), 427-439 \ (2009).

Hanson, A. J.: Visualizing quaternions, Elsevier Morgan Kaufmann Publishers,
(2006).

Jia, Y.B.: Quaternion and Rotation, Com S 477/577 Notes, (2017).

Jouget, P.: S\'{e}curit\'{e} et performance de dispositifs de distribution
quantique de cl\'{e}s \`{a} variables continues, PhD Thesis, TELECOM
ParisTech, (2013).

Klco, P., Smetana, M., Kollarik, M., Tatar, M.: Application of Octonions in
the Cough Sounds Classification , Advances in Applied Science Research, 
\textbf{8}(2), 30-37 (2017).

Kostrikin, A. I., Shafarevich, I.R. (Eds): Algebra VI, Springer-Verlag,
(1995).

Li, X.M.: Hyper-Complex Numbers and its Applications in Digital Image
Processing, Seminars and Distinguished Lectures, (2011).

Schafer, R. D.: An Introduction to Nonassociative Algebras\textit{,}
Academic Press, New-York, (1966).

Snopek, K., M.: Quaternions and Octonions in Signal Processing -
Fundamentals and Some New Results, Przeglad Telekomunikacyjny - Wiadomo\'{s}%
ci Telekomunikacyjne, SIGMA NOT, \textbf{134}(6), 619-622 \ (2015).

Unger, T., Markin, N.: Quadratic Forms and Space-Time Block Codes from
Generalized Quaternion and Biquaternion Algebras, Information Theory, IEEE
Transactions, \textbf{57(}9\textbf{)}, 6148-6156 \ (2011).%
\begin{equation*}
\end{equation*}

Cristina FLAUT

{\small Faculty of Mathematics and Computer Science, Ovidius University,}

{\small Bd. Mamaia 124, 900527, CONSTANTA, ROMANIA}

{\small http://cristinaflaut.wikispaces.com/;
http://www.univ-ovidius.ro/math/}

{\small e-mail: cflaut@univ-ovidius.ro; cristina\_flaut@yahoo.com}

\end{document}